\documentclass[a4paper]{article}

\usepackage{graphicx}

\usepackage{newtxtext} 
\usepackage[varvw]{newtxmath}

\usepackage{tikz}
\usetikzlibrary{decorations.pathreplacing}
\usetikzlibrary{arrows.meta}

\usepackage{xifthen}
\newcommand{\pateraBlock}[7]{

    \ifthenelse{#2 = 1}{\draw (0, -0.25*\s) -- (1*\s, -0.25*\s);}{\ifthenelse{#2 = 2}{\draw[dashed] (0, -0.27*\s) -- (1*\s, -0.27*\s);}{}}
    \draw (0, -0.25*\s) -- (0, 0);
    \draw (1*\s, -0.25*\s) -- (1*\s, 0);
    
    \ifthenelse{#3 = 1}{\draw (-0.25*\s, 0) -- (-0.25*\s, 1*\s);}{\ifthenelse{#3 = 2}{\draw[dashed] (-0.27*\s, 0) -- (-0.27*\s, 1*\s);}{}}
    \draw (-0.25*\s, 0) -- (0, 0);
    \draw (-0.25*\s, 1*\s) -- (0, 1*\s);

    \ifthenelse{#4 = 1}{\draw (1.25*\s, 0) -- (1.25*\s, 1*\s);}{\ifthenelse{#4 = 2}{\draw[dashed] (1.27*\s, 0) -- (1.27*\s, 1*\s);}{}}
    \draw (1.25*\s, 0) -- (1*\s, 0);
    \draw (1.25*\s, 1*\s) -- (1*\s, 1*\s);

    \ifthenelse{#5 = 1}{\draw (0, 1.25*\s) -- (1*\s, 1.25*\s);}{\ifthenelse{#5 = 2}{\draw[dashed] (0, 1.27*\s) -- (1*\s, 1.27*\s);}{}}
    \draw (0, 1*\s) -- (0, 1.25*\s);
    \draw (1*\s, 1*\s) -- (1*\s, 1.25*\s);
    
    \draw[red!10, fill = red!10] (0, 0) rectangle (1*\s, 1*\s);
    
    \ifthenelse{#6 = 1}{\node at (0.5*\s, 0.65*\s) {$\mu_{#1}$};}{\ifthenelse{#6 = 2}{\node at (0.5*\s, 0.65*\s) {$\hmu$};}}

    \ifthenelse{#7 = 1}{\node at (0.5*\s, 0.35*\s) {$\Omega^b_{#1}$};}{\ifthenelse{#7 = 2}{\node at (0.5*\s, 0.35*\s) {$\hOmega^b_{#1}$};}}

}

\newcommand{\pateraWing}[4]{

    \ifthenelse{#1 = 1}{\draw[gray!5, fill = gray!5] (0, -0.25*\s) rectangle (1*\s, 0);}{\ifthenelse{#1 = 2}{\draw[blue!5, fill = blue!5] (0, -0.25*\s) rectangle (1*\s, 0);}}{\draw[gray!5, fill = gray!5] (0, -0.20*\s) rectangle (1*\s, 0);} 
    \ifthenelse{#2 = 1}{\draw[gray!5, fill = gray!5] (0, 1*\s) rectangle (1*\s, 1.25*\s);}{\ifthenelse{#2 = 2}{\draw[blue!5, fill = blue!5] (0, 1*\s) rectangle (1*\s, 1.25*\s);}}{\draw[gray!5, fill = gray!5] (0, 1*\s) rectangle (1*\s, 1.20*\s);} 
    \ifthenelse{#3 = 1}{\draw[gray!5, fill = gray!5] (-0.25*\s, 0) rectangle (0, 1*\s);}{\ifthenelse{#3 = 2}{\draw[yellow!5, fill = yellow!5] (-0.25*\s, 0) rectangle (0, 1*\s);}}{\draw[gray!5, fill = gray!5] (-0.20*\s, 0) rectangle (0, 1*\s);} 
    \ifthenelse{#4 = 1}{\draw[gray!5, fill = gray!5] (1*\s, 0) rectangle (1.25*\s, 1*\s);}{\ifthenelse{#4 = 2}{\draw[yellow!5, fill = yellow!5] (1*\s, 0) rectangle (1.25*\s, 1*\s);}}{\draw[gray!5, fill = gray!5] (1*\s, 0) rectangle (1.20*\s, 1*\s);} 

}

\newcommand{\bm}[1]{\text{\boldmath $#1$\unboldmath}}
\newcommand{\bx}{\bm{x}}
\newcommand{\bn}{\bm{n}}
\newcommand{\bmu}{\boldsymbol{\mu}}
\newcommand{\pgd}{{\text{\tiny{PGD}}}}
\newcommand{\bLambda}{\boldsymbol{\Lambda}}

\usepackage{a4wide}
\begin{document}

\begin{center}

\begin{Large}
\textbf{Local surrogate models with reduced dimensionality via overlapping domain decomposition and proper generalized decomposition}
\end{Large}

\medskip

Marco Discacciati$^1$, Ben J. Evans$^1$, Matteo Giacomini$^{2,3}$

\medskip

${}^1$ Department of Mathematical Sciences, Loughborough University, Epinal Way, Loughborough, UK, m.discacciati@lboro.ac.uk, b.j.evans@lboro.ac.uk.

${}^2$ Laboratori de C\`alcul Numeric (LaC\`aN), E.T.S. de Ingenier\'ia de Caminos, Canales y Puertos, Universitat Polit\`ecnica de
Catalunya, Barcelona, Spain

${}^3$ Centre Internacional de M\`etodes Num\`erics en Enginyeria (CIMNE), Barcelona, Spain, matteo.giacomini@upc.edu.

\end{center}

\bigskip

\centerline{\textbf{Abstract}}

{We propose an efficient algorithm that combines overlapping domain decomposition and proper generalized decomposition (PGD) to construct surrogate models of linear elliptic parametric problems. The technique is composed of an offline and an online phase that can be implemented in a fully non-intrusive way. The online phase relies on a substructured algebraic formulation of the alternating Schwarz method, while the offline phase exploits the linearity of the boundary value problem to characterize a \textit{PGD basis} and generate local surrogate models, with minimal parametric dimensionality, in each subdomain. Numerical results show the efficiency of the proposed methodology.}

\section{Introduction}

The combination of domain decomposition (DD) methods and model order reduction (MOR) techniques can lead to efficient algorithms to solve parametric partial differential equations for challenging applications such as, e.g., the construction of digital twins involving multi-physics models. Indeed, DD can split the target problem into local subproblems characterized by a single physics and by a reduced number of parameters and degrees of freedom that can be more easily handled by MOR techniques (see, e.g., \cite{buhr-2021,deg-stokesDarcy-2024}).

In \cite{deg-cmame-2024}, we proposed a computational framework based on a substructured overlapping DD algorithm to couple local surrogate models generated by proper generalized decomposition (PGD). This is a physics-based a priori MOR method (see, e.g., \cite{ckl-pgdBook}) that does not require any sampling of the parametric space and that can be very conveniently implemented in a completely non-intrusive way \cite{epgd-2020}. By exploiting the linearity of the underlying parametric partial differential equation, in the offline phase, we construct local surrogate PGD models that incorporate arbitrary traces of the solution at the subdomain interfaces, to be used in the online phase to perform alternating Schwarz iterations.
However, the construction of such local surrogate models requires to suitably cluster the interface nodes, to appropriately define the parametric values of the traces, and to solve parametric subproblems with combinations of boundary conditions that could be not physically relevant.

To overcome these drawbacks, in this work, we improve the computational framework  of \cite{deg-cmame-2024} by reducing the parametric dimension of the offline local problems. This avoids the burden of the implementation of the clustering procedure, and it allows to compute local \textit{PGD basis} functions that can be easily linearly combined during the online phase to obtain the local surrogate models with the interface trace needed by the Schwarz algorithm. This novel approach simplifies the overall implementation of the DD-PGD procedure and it results in significantly lower computational costs both offline and online.

After introducing the model parametric problem in Sect.~\ref{sec:setting}, we present the novel approach to construct the local surrogate models in Sect.~\ref{sec:offline}, and we show how to perform the online phase in Sect.~\ref{sec:onlinePhase}. Finally, numerical results in Sect.~\ref{sec:numericalResults} assess the performance of the algorithm compared to its counterpart in \cite{deg-cmame-2024}.

\section{Problem formulation and PGD-accelerated Schwarz method}
\label{sec:setting}

In this section, we formulate the alternating Schwarz method for a parametric model problem, and we explain how the algorithm can be split into an offline and an online phase to exploit the benefits of MOR. For the sake of simplicity, we consider a parametric Poisson problem with homogeneous Dirichlet boundary conditions, we assume that the domain is parameter-independent, and we introduce a domain decomposition into two overlapping subdomains only. The extension to more general linear elliptic operators with arbitrary boundary conditions and decompositions into multiple subdomains (without cross-points) can be found in \cite{deg-cmame-2024}.

Let $\Omega \subset \mathbb{R}^d$ ($d = 1, 2, 3$) be an open bounded domain with Lipschitz boundary $\partial\Omega$, and let $\bmu =(\mu^1,\ldots,\mu^P) \in \mathcal{P}$ be a tuple of $P \in \mathbb{N}$ problem parameters with $\mathcal{P} = \mathcal{I}^1 \times \dots \times \mathcal{I}^P \subset \mathbb{R}^P$ with $\mathcal{I}^p$ compact ($p=1,\ldots,P$). Moreover, for all $\bmu \in \mathcal{P}$, let $\nu(\bmu) >0$ be a positive parametric diffusion coefficient, and $s(\bmu) \in L^2(\Omega)$ a parametric source term. Then, consider the well-posed linear elliptic parametric boundary value problem: for all $\bmu \in \mathcal{P}$, find $u(\bmu)$ such that 
\begin{equation}
	\label{eq:globalProb}
	\begin{array}{rcll}
		-\nabla\cdot(\nu(\bmu)\nabla u(\bmu)) &=& s(\bmu) &\quad \text{in } \Omega,\\
		u(\bmu) &=& 0 & \quad \text{on } \partial\Omega.
	\end{array}
\end{equation}
The domain $\Omega$ is split into two overlapping subdomains $\Omega_i \subset \Omega$ ($i = 1, 2$) such that $\Omega_1 \cup \Omega_2 = \Omega$ and $\Omega_1 \cap \Omega_2 
\neq \emptyset$, and, for $i=1,2$, let $\Gamma_i = \partial\Omega_i \setminus \partial\Omega$. Considering this decomposition, problem~\eqref{eq:globalProb} can be equivalently reformulated in the multi-domain form:
for all $\bmu \in \mathcal{P}$, find $u_i(\bmu)$ ($i = 1, 2$) such that
\begin{equation}
	\label{eq:multiDom}
	\begin{array}{rcll}
		-\nabla\cdot(\nu(\bmu)\nabla u_i(\bmu)) &=& s_i(\bmu) & \quad \text{in } \Omega_i,\\
		u_i(\bmu) &=& 0 & \quad \text{on } \partial\Omega_i \cap \partial\Omega,\\
		u_1(\bmu) &=& u_2(\bmu) & \quad \text{on } \Gamma_1 \cup \Gamma_2,
	\end{array}
\end{equation}
where $s_i(\bmu)$ denotes the restriction of $s(\bmu)$ to $\Omega_i$, and \eqref{eq:multiDom}$_3$ ensures the continuity of the local solutions $u_1(\bmu)$ and $u_2(\bmu)$ across the interfaces $\Gamma_1$ and $\Gamma_2$.

\subsection{Parametric overlapping Schwarz method}\label{sec:parametricAS}

Using the multi-domain formulation~\eqref{eq:multiDom}, we can adapt the classical overlapping alternating Schwarz method to the case of a parametric problem. More precisely, let $\bar{\bmu} \in \mathcal{P}$ be a fixed set of parameters, and let $\lambda_1^{(0)}$ be a suitable trace on the interface $\Gamma_1$. Then, for $k \geq 1$ until convergence when $\|u_1^{(k)}\!(\bar{\bmu}) - u_2^{(k)}\!(\bar{\bmu})\|_{\Gamma} < tol$, for a suitable tolerance $tol$ and a suitable norm $\|\cdot\|_\Gamma$ on $\Gamma_1\cup\Gamma_2$,
\begin{itemize}
 \item 
    Find $u_1^{(k)}\!(\bar{\bmu})$ such that 
	 \begin{equation}\label{eq:localProblem1}
		\begin{array}{rcll}
			-\nabla \cdot (\nu(\bar{\bmu}) \nabla u_1^{(k)}(\bar{\bmu})) &=& s_1(\bar{\bmu}) & \quad \text{in } \Omega_1,\\
			u_1^{(k)}\!(\bar{\bmu}) &=& 0 & \quad \text{on } \partial\Omega_1 \cap \partial \Omega,\\
			u_1^{(k)}\!(\bar{\bmu}) &=& \lambda_1^{(k-1)} & \quad \text{on } \Gamma_1 .
		\end{array}
	\end{equation}
	\label{alg:firstStep}
 \item
 Set $\lambda_2^{(k)} = u_1^{(k)}\!(\bar{\bmu})|_{\Gamma_2}$.
 \item
 Find $u_2^{(k)}\!(\bar{\bmu})$ such that 
	\begin{equation}\label{eq:localProblem2}
		\begin{array}{rcll}
			-\nabla \cdot (\nu(\bar{\bmu}) \nabla u_2^{(k)}(\bar{\bmu})) &=& s_2(\bar{\bmu}) & \quad \text{in } \Omega_2,\\
			u_2^{(k)}\!(\bar{\bmu}) &=& 0 & \quad \text{on } \partial\Omega_2 \cap \partial \Omega,\\
			u_2^{(k)}\!(\bar{\bmu}) &=& \lambda_2^{(k)} & \quad \text{on } \Gamma_2.
		\end{array}
	\end{equation}
 \item
 Set $\lambda_1^{(k)}= u_2^{(k)}\!(\bar{\bmu})|_{\Gamma_1}$.
\end{itemize}

To efficiently obtain the solution of the parametric problem~\eqref{eq:globalProb} for any set of parameters $\bmu \in \mathcal{P}$, the execution of the alternating Schwarz algorithm can be split into an offline phase and an online phase. In the former, PGD is used to precompute local surrogate models that can handle arbitrary Dirichlet boundary conditions at the interfaces, thus reducing the computational cost of the online Schwarz iterations.

More precisely, in the \textit{offline} phase, the local parametric problems~\eqref{eq:localProblem1} and~\eqref{eq:localProblem2} are solved by PGD so that, at the end of this phase, local surrogate models $u^\pgd_i$ are available. The local surrogate models $u^\pgd_i$ feature arbitrary traces $\lambda_i$ at the interfaces, and this is achieved by parametrizing them through a set of auxiliary parameters, say $\bLambda_i$. As we will see in Sect.~\ref{sec:offline}, the linearity of the local problems is exploited to avoid incurring the curse of dimensionality after introducing such additional parameters.

In the \textit{online} phase, the alternating Schwarz iterations are performed for a fixed set of parametric values $\bar{\bmu} \in \mathcal{P}$, and the solution of the local problems~\eqref{eq:localProblem1} and~\eqref{eq:localProblem2} is replaced by the \textit{evaluation} of the precomputed local surrogate models $u^\pgd_i$ at specific instances of the interface parameters $\bLambda_i$. Since no problems are solved at this stage, the Schwarz algorithm can be executed in real time. Details are provided in Sect.~\ref{sec:onlinePhase}.

\section{Offline phase: Computation of the local surrogate models}
\label{sec:offline}

In the offline phase of the algorithm, PGD is used to solve the following local parametric problem: for all $\bmu \in \mathcal{P}$, find $u_i(\bmu)$ $(i=1,2)$ such that
\begin{equation}
	\label{eq:localProb}
	\begin{array}{rcll}
		- \nabla \cdot ( \nu(\bmu) \nabla u_{i}(\bmu) ) &=& s_i(\bmu) & \quad \text{in } \Omega_i,\\
		u_i(\bmu) &=& 0 & \quad \text{on } \partial \Omega_i \cap \partial \Omega,\\
		u_i(\bmu) &=& \lambda_i & \quad \text{on } \Gamma_i,
	\end{array}
\end{equation}
where $\lambda_i$ is a space-dependent function that represents an arbitrary Dirichlet boundary condition at the interface $\Gamma_i$. A suitable parametrization is introduced to handle the arbitrariness of the trace $\lambda_i$. More precisely, considering that \eqref{eq:localProb} is solved, e.g., by the finite element method, we express $\lambda_i$ as a linear combination of suitable finite element basis functions on $\Gamma_i$:
\begin{equation}\label{eq:lambdaRep}
 \lambda_i(\bx) = \sum_{q = 1}^{N_{\Gamma_i}} \Lambda^q_i \, \eta^q_i(\bx) \, .
\end{equation}
Here, $N_{\Gamma_i}$ represents the dimension of the discrete trace space generated by the basis functions $\eta^q_i(\bx)$ on the interface $\Gamma_i$, while $\bLambda_i = (\Lambda^1_i, \ldots, \Lambda^{N_{\Gamma_i}}_i) \in \mathbb{R}^{N_{\Gamma_i}}$ are the coefficients of the linear combination. Following \cite{deg-cmame-2024}, we can choose $\eta^q_i(\bx)$ to be the non-null restriction on $\Gamma_i$ of the standard piecewise Lagrange polynomial basis functions, say $\varphi^q_i(\bx)$, that discretize \eqref{eq:localProb} in $\Omega_i$, that is, $\eta^q_i(\bx) = \varphi^q_i(\bx)|_{\Gamma_i}$. Therefore, the basis functions $\eta^q_i(\bx)$ are linearly independent, they provide a partition of unity on $\Gamma_i$, and they satisfy $\eta^q_i(\bx_i^m) = \delta_{qm}$ for all nodes $\bx_i^m$ on the interface $\Gamma_i$, $\delta_{qm}$ being the Kronecker delta. (Although this choice of the basis functions on $\Gamma_i$ is very convenient for implementation purposes, other choices are possible and they can be accommodated in the overall structure of the algorithm.)

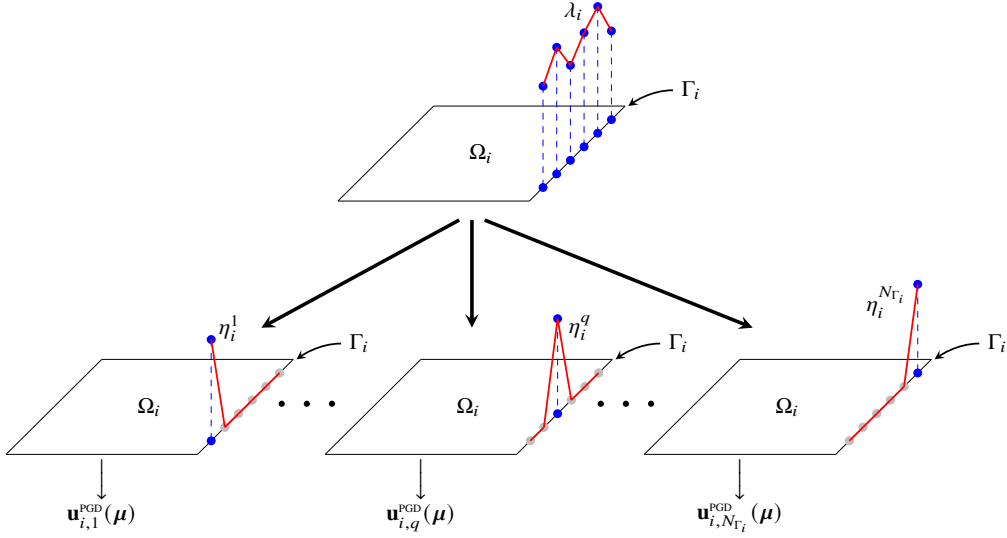
\begin{figure}[!ht]
    \centering
    \resizebox{0.9\textwidth}{!}{
    \begin{tikzpicture}
			\newcommand{\intNodes}[3]{
				\draw[dashed, #3] (\longSide + #1 * \delta / 7, \ang * #1 * \delta / 7 + #2) -- (\longSide + #1 * \delta / 7, #1 * \ang * \delta / 7 );
				\fill[#3] (\longSide + #1 * \delta / 7, \ang * #1 * \delta / 7 + #2) circle (0.07 cm);
				\fill[#3] (\longSide + #1 * \delta / 7, #1 * \ang * \delta / 7 ) circle (0.07 cm);
			}
			
			\newcommand{\intFun}[3]{
				\draw[thick, red] (\longSide + #1 * \delta / 7, #1 * \ang * \delta / 7 + #2 ) -- (\longSide + #1 * \delta / 7 + \delta / 7, #1 * \ang * \delta / 7 +  \ang * \delta / 7 + #3);
			}
			
			\def\ang{1}
			\def\delta{1.5}
			\def\longSide{3}
			
			\begin{scope}[shift={(0.2, 0)}]
				\draw (0, 0) -- (\delta, \delta*\ang);
				\draw (0, 0) -- (\longSide, 0);
				\draw(\delta, \delta*\ang) -- (\delta + \longSide, \delta*\ang);
				\draw(\longSide, 0) -- (\delta + \longSide, \delta*\ang);
				
				\node at (\longSide / 2 + \delta / 2, \delta * \ang / 2) {$\Omega_i$};
				
				\draw[-stealth, thick] (\longSide + \delta + 0.75, \delta * \ang + 0.25) node[anchor=west] {$\Gamma_i$} parabola (\longSide + \delta + 0.05, \delta * \ang);
				
				\intNodes{1}{1.6}{blue};
				\intNodes{2}{2}{blue};
				\intNodes{3}{1.5}{blue};	
				\intNodes{4}{1.8}{blue};
				\intNodes{5}{2}{blue};
				\intNodes{6}{1.4}{blue};
				
				\intFun{1}{1.6}{2};
				\intFun{2}{2}{1.5};
				\intFun{3}{1.5}{1.8};
				\intFun{4}{1.8}{2};
				\intFun{5}{2}{1.4};
			\end{scope}
			
			\begin{scope}[shift={(-5, -4)}]
				\draw (0, 0) -- (\delta, \delta*\ang);
				\draw (0, 0) -- (\longSide, 0);
				\draw(\delta, \delta*\ang) -- (\delta + \longSide, \delta*\ang);
				\draw(\longSide, 0) -- (\delta + \longSide, \delta*\ang);
				
				\node at (\longSide / 2 + \delta / 2, \delta * \ang / 2) {$\Omega_i$};
				
				\draw[-stealth, thick] (\longSide + \delta + 0.75, \delta * \ang + 0.25) node[anchor=west] {$\Gamma_i$} parabola (\longSide + \delta + 0.05, \delta * \ang);
				
				\intNodes{1}{1.6}{blue};
				\intNodes{2}{0}{lightgray};
				\intNodes{3}{0}{lightgray};
				\intNodes{4}{0}{lightgray};
				\intNodes{5}{0}{lightgray};
				\intNodes{6}{0}{lightgray};

				\intFun{1}{1.6}{0};
				\intFun{2}{0}{0};
				\intFun{3}{0}{0};
				\intFun{4}{0}{0};
				\intFun{5}{0}{0};
			\end{scope}
		
			\begin{scope}[shift={(0, -4)}]
				\draw (0, 0) -- (\delta, \delta*\ang);
				\draw (0, 0) -- (\longSide, 0);
				\draw(\delta, \delta*\ang) -- (\delta + \longSide, \delta*\ang);
				\draw(\longSide, 0) -- (\delta + \longSide, \delta*\ang);
				
				\node at (\longSide / 2 + \delta / 2, \delta * \ang / 2) {$\Omega_i$};
				
				\draw[-stealth, thick] (\longSide + \delta + 0.75, \delta * \ang + 0.25) node[anchor=west] {$\Gamma_i$} parabola (\longSide + \delta + 0.05, \delta * \ang);
				
				\intNodes{1}{0}{lightgray};
				\intNodes{2}{0}{lightgray};
				
				\intNodes{3}{1.5}{blue};
				
				\intNodes{4}{0}{lightgray};
				\intNodes{5}{0}{lightgray};
				\intNodes{6}{0}{lightgray};

				\intFun{1}{0}{0};
				\intFun{2}{0}{1.5};
				\intFun{3}{1.5}{0};
				\intFun{4}{0}{0};
				\intFun{5}{0}{0};
			\end{scope}
		
			\begin{scope}[shift={(5, -4)}]
				\draw (0, 0) -- (\delta, \delta*\ang);
				\draw (0, 0) -- (\longSide, 0);
				\draw(\delta, \delta*\ang) -- (\delta + \longSide, \delta*\ang);
				\draw(\longSide, 0) -- (\delta + \longSide, \delta*\ang);
				
				\node at (\longSide / 2 + \delta / 2, \delta * \ang / 2) {$\Omega_i$};
				
				\draw[-stealth, thick] (\longSide + \delta + 0.75, \delta * \ang + 0.25) node[anchor=west] {$\Gamma_i$} parabola (\longSide + \delta + 0.05, \delta * \ang);
				
				\intNodes{1}{0}{lightgray};
				\intNodes{2}{0}{lightgray};	
				\intNodes{3}{0}{lightgray};
				\intNodes{4}{0}{lightgray};
				\intNodes{5}{0}{lightgray};
				\intNodes{6}{1.4}{blue};				
				
				\intFun{1}{0}{0};
				\intFun{2}{0}{0};
				\intFun{3}{0}{0};
				\intFun{4}{0}{0};
				\intFun{5}{0}{1.4};
			\end{scope}
			
			\draw[- stealth, ultra thick] (2.1, -0.3) -- (-1 , -2);
			\draw[- stealth, ultra thick] (2.3, -0.3) -- (2.3 , -2);
			\draw[- stealth, ultra thick] (2.5, -0.3) -- (6.7 , -2);
			
			\node at (-0.2, -3.2) {\Huge \ldots};
			\node at (4.8, -3.2) {\Huge \ldots};
			
			\node at (-3.5, -4.4) {$\Big\downarrow$};
			\node at ( 1.5, -4.4) {$\Big\downarrow$};
			\node at ( 6.5, -4.4) {$\Big\downarrow$};
			
			\node at (-3.5, -5.0) {$\mathbf{u}_{i,1}^{\pgd}(\bmu)$};
			\node at ( 1.5, -5.0) {$\mathbf{u}_{i,q}^{\pgd}(\bmu)$};
			\node at ( 6.5, -5.0) {$\mathbf{u}_{i,N_{\Gamma_i}}^{\pgd}(\bmu)$};

                \node at ( 3.9,  3.0) {$\lambda_i$};
                \node at (-1.5, -2.0) {$\eta_i^1$};
                \node at ( 4.0, -2.0) {$\eta_i^q$};
                \node at ( 8.8, -1.6) {$\eta_i^{N_{\Gamma_i}}$};
			
		\end{tikzpicture}
  }
    \caption{Partition of the interface nodes as a collection of single independent boundary parameters.}
    \label{fig:activeBdryNodes}
\end{figure}

The coefficients $\Lambda^q_i$ could be considered as additional parameters of the local problem \eqref{eq:localProb}, and the latter be solved by PGD. However, in general we expect $N_{\Gamma_i} \gg 1$, and it is well-known that a limitation of MOR methods is that they cannot efficiently handle a large number of parameters. To overcome this difficulty, we exploit the linearity of the operator that governs the boundary value problem and we introduce the following $N_{\Gamma_i}+1$ problems:
\begin{itemize}
\item
for all $\bmu \in \mathcal{P}$, find $u_{i,0} (\bmu)$ such that
\begin{subequations}\label{eq:subProb}
\begin{equation}\label{eq:sourceProb}
\begin{array}{rcll}
- \nabla \cdot ( \nu(\bmu) \nabla {u}_{i, 0}(\bmu)) &=& s_i(\bmu) & \quad \text{in } \Omega_i,\\
u_{i, 0}(\bmu) &=& 0 & \quad \text{on } \partial\Omega_i \cap \partial\Omega,\\
u_{i, 0}(\bmu) &=& 0 & \quad \text{on } \Gamma_i ;
\end{array}
\end{equation}
\item
for $q=1,\ldots,N_{\Gamma_i}$ and for all $\bmu \in \mathcal{P}$, find $u_{i, q}(\bmu)$ such that 
\begin{equation}\label{eq:boundaryProbs}
\begin{array}{rcll}
- \nabla \cdot ( \nu(\bmu) \nabla u_{i, q}(\bmu)) &=& 0 &\quad \text{in } \Omega_i,\\
u_{i, q}(\bmu) &=& 0 &\quad \text{on } \partial\Omega_i \cap \partial\Omega,\\
u_{i, q}(\bmu) &=& \eta^q_i & \quad \text{on } \Gamma_i .
\end{array}
\end{equation}
\end{subequations}
\end{itemize}
A schematic representation of this idea is presented in Figure \ref{fig:activeBdryNodes}.
Problem \eqref{eq:sourceProb} and all problems \eqref{eq:boundaryProbs} can be solved independently, and, being associated with the linearly independent trace functions $\eta_i^q$, the solutions $\{u_{i,q}(\bmu)\}_{q=1,\ldots,N_{\Gamma_i}}$ of \eqref{eq:boundaryProbs} can be considered as a \textit{PGD basis} that can be used to construct the local surrogate model in $\Omega_i$. More precisely, for all $\bmu \in \mathcal{P}$ and for a given trace function $\lambda_i$ on $\Gamma_i$, upon considering the expansion \eqref{eq:lambdaRep}, the solution of \eqref{eq:localProb} can be obtained as
\begin{equation}\label{eq:solutionSplit}
u_i(\bmu, \bLambda_i) = u_{i,0} (\bmu) + \sum_{q=1}^{N_{\Gamma_i}} \Lambda_i^q \, u_{i,q}(\bmu) .
\end{equation}

This approach is different from the one proposed in \cite{deg-cmame-2024} and it presents several computational advantages, both in the offline and in the online phase. We focus here on the offline phase and we will discuss the online phase in Sect.~\ref{sec:onlinePhase}. In \cite{deg-cmame-2024}, to overcome the difficulty of handling too many interface parameters $\Lambda_i^q$, these were clustered in $N_i$ sufficiently small disjoint sets $\mathcal{N}_i^j$ ($j=1,\ldots,N_i$) of \textit{active boundary parameters} with $\bigcup_{j=1,\ldots,N_i} \mathcal{N}_i^j = \{1,\ldots,N_{\Gamma_i}\}$ and with $\text{card}(\mathcal{N}_i^j) \ll N_{\Gamma_i}$. Then, \eqref{eq:lambdaRep} was be replaced by 
\begin{equation}\label{eq:splittingBoundaryParameters}
\lambda_i = 
\sum_{q \in \mathcal{N}^1_i} \Lambda^q_i \, \eta^q_i + 
\sum_{q \in \mathcal{N}^2_i} \Lambda^q_i \, \eta^q_i + \ldots +
\sum_{q \in \mathcal{N}^{N_i}_i} \Lambda^q_i \, \eta^q_i ,
\end{equation}
and, instead of \eqref{eq:boundaryProbs}, one would solve the $N_i$ independent problems: for $j=1,\ldots,N_i$, for all $\bmu \in \mathcal{P}$ and for all $\bLambda_i^j = (\Lambda_i^q)_{q\in \mathcal{N}_i^j} \in \mathcal{Q}_i^j$, find 
\begin{equation}\label{eq:boundaryProbsCluster}
\begin{array}{rcll}
- \nabla \cdot ( \nu(\bmu) \nabla u_{i, j}(\bmu,\bLambda_i^j)) &=& 0 &\quad \text{in } \Omega_i,\\
u_{i, j}(\bmu,\bLambda_i^j) &=& 0 &\quad \text{on } \partial\Omega_i \cap \partial\Omega,\\
u_{i, j}(\bmu,\bLambda_i^j) &=& \sum_{q \in \mathcal{N}_i^j} \Lambda_i^q \eta^q_i & \quad \text{on } \Gamma_i .
\end{array}
\end{equation}
Here, $\mathcal{Q}_i^j = \bigtimes_{q \in \mathcal{N}_i^j} \mathcal{J}_i^q$ is the space of boundary parameters with each $\mathcal{J}_i^q \subset \mathbb{R}$ ($q = 1, \dots , N_{\Gamma_i}$) a compact set. Finally, the solution of \eqref{eq:localProb} is obtained as
\begin{equation}\label{eq:solutionSplitCluster}
u_i(\bmu, \bLambda_i) = u_{i,0} (\bmu) + \sum_{j=1}^{N_i} u_{i,j}(\bmu, \bLambda_i^j) .
\end{equation}
In the case of clustered interface nodes, one must solve a smaller number $N_i \ll N_{\Gamma_i}$ of local problems \eqref{eq:boundaryProbsCluster}, each depending both on $\bmu$ and on the interface parameters $\bLambda_i^j$, while \eqref{eq:boundaryProbs} only involves the problem parameter $\bmu$. Hence, the proposed approach leads to local problems with significantly reduced dimensionality, at the expense of computing $N_{\Gamma_i}$ local surrogate models instead of $N_i$. Nonetheless, it is worth recalling that the local problems \eqref{eq:boundaryProbs} are independent from one another and can be easily solved in parallel. Moreover, before solving \eqref{eq:boundaryProbsCluster}, a suitable interface parametric space $\mathcal{Q}_i^j$ must be identified to represent an arbitrary trace that is non-zero at the nodes corresponding to the indices in $\mathcal{N}_i^j$. In principle, each intervals $\mathcal{J}_i^q$ depends on $\bmu$ and it should be defined considering physical information about the solution of problem \eqref{eq:globalProb} (e.g., the maximum and minimum values that $u(\bmu)$ can attain within $\Omega$). In the approach proposed in this work, this is not needed, as the auxiliary problems \eqref{eq:boundaryProbs} are independent of $\bLambda_i^j$. In particular, it should be noted that not all possible combinations of the parameters $\Lambda_i^q$ considered in \eqref{eq:boundaryProbsCluster} may actually be significant for the problem at hand. For instance, considering very different values of $\Lambda_i^q$ at adjacent nodes on $\Gamma_i$ can represent a highly oscillatory trace function, that it is unlikely to be relevant for the solution of an elliptic problem in a smooth domain with sufficiently regular data. Moreover, solving boundary value problems with such localised features can be challenging for PGD due to the possibly large number of modes needed to correctly represent the behaviour of the solution. Therefore, the offline computational cost associated with the clustering procedure can become unnecessarily high, as numerical simulations will show in Sect.~\ref{sec:numericalResults}. Finally, no clustering procedure is needed when working with \eqref{eq:sourceProb} and \eqref{eq:boundaryProbs}, thus making the implementation simpler.

\subsection{Construction of the local surrogate models using PGD}\label{sec:surrogateModels}

In this section, we outline the procedure to compute the solutions of the local problems \eqref{eq:sourceProb} and \eqref{eq:boundaryProbs} using PGD, and we refer to \cite{deg-cmame-2024} for more details.

Following the standard approach in PGD \cite{ckl-pgdBook}, we assume that all data are given in separated form as the sum of products of functions that depend either on the spatial coordinate $\bx$ or on the parameters $\bmu$. The contribution of the Dirichlet boundary condition in \eqref{eq:boundaryProbs} is handled by introducing ad-hoc, sufficiently smooth modes as usually done in the PGD context. Then, the solution of \eqref{eq:sourceProb} and the solution, say $v_{i,q}(\bmu)$, of \eqref{eq:boundaryProbs}, after accounting for the Dirichlet boundary condition on $\Gamma_i$, are written in the PGD expansion
\begin{equation*}
\resizebox{0.95\hsize}{!}{$\displaystyle
    u_{i,0}(\bmu) \approx u_{i,0}^{\pgd}(\bmu) = \sum_{m=1}^{M_0} V_{i,0}^m(\bx) \phi_{i,0}^m (\bmu) \, ,
    \quad  
    v_{i,q}(\bmu) \approx v_{i,q}^{\pgd}(\bmu) = \sum_{m=1}^{M_q} V_{i,q}^m(\bx) \phi_{i,q}^m (\bmu) \, ,
$}
\end{equation*}
with $V_{i,0}^m$ and $V_{i,j}^m$ being the $m$-th spatial modes, and $\phi_{i,0}^m (\bmu)$ and $\phi_{i,q}^m (\bmu)$ the parametric modes. The spatial modes are discretized using continuous Galerkin finite elements, while pointwise collocation is used to approximate the parametric modes. The number of modes $M_0$ and $M_q$ is automatically determined by a greedy procedure. Since all problems \eqref{eq:sourceProb} and \eqref{eq:boundaryProbs} are independent, we compute their PGD approximations separately using the encapsulated PGD library \cite{epgd-2020}. At the end of this procedure, we obtain the discrete PGD approximations $\mathbf{u}_{i,0}^\pgd(\bmu)$ and $\mathbf{u}_{i,q}^\pgd(\bmu)$ of the solutions of problems \eqref{eq:sourceProb} and \eqref{eq:boundaryProbs}, respectively.

\section{Online phase: PGD-accelerated Schwarz iterations}
\label{sec:onlinePhase}

In this section, we present an efficient algorithm to construct the global solution of \eqref{eq:globalProb} for fixed values $\bar{\bmu} \in \mathcal{P}$ of the physical parameter.

First of all, considering a suitable conforming finite element approximations of the local problems \eqref{eq:localProblem1} and \eqref{eq:localProblem2}, we recall (see, e.g., \cite{smith-ddBook}) that one step of the overlapping Schwarz method outlined in Sect.~\ref{sec:parametricAS} corresponds to one Gauss-Seidel iteration to solve the linear system
\begin{equation}\label{eq:systemSchwarz}
\begin{pmatrix}
{A}_{\Omega_1} & {A}_{\Gamma_1} & {0} & {0} \\
{0} & {I}_{\Gamma_1} & - \mathcal{I}_{\Omega_2\to\Gamma_1} & {0} \\
{0} & {0} & {A}_{\Omega_2} & {A}_{\Gamma_2} \\
- \mathcal{I}_{\Omega_1\to\Gamma_2} & {0} & {0} & {I}_{\Gamma_2}
\end{pmatrix}
\begin{pmatrix}
\mathbf{u}_{\Omega_1}\!(\bar{\bmu}) \\
\mathbf{u}_{\Gamma_1}\!(\bar{\bmu}) \\
\mathbf{u}_{\Omega_2}\!(\bar{\bmu}) \\
\mathbf{u}_{\Gamma_2}\!(\bar{\bmu})
\end{pmatrix}
=
\begin{pmatrix}
\mathbf{f}_{\Omega_1}\!(\bar{\bmu}) \\
\mathbf{0} \\
\mathbf{f}_{\Omega_2}\!(\bar{\bmu}) \\
\mathbf{0}
\end{pmatrix}\,.
\end{equation}
Here, $\mathbf{u}_{\Omega_i}(\bar{\bmu})$ and $\mathbf{u}_{\Gamma_i}(\bar{\bmu})$ denote the vectors of nodal values of the unknown solutions inside the domain $\Omega_i$ and on $\Gamma_i$, respectively, whereas ${I}_{\Gamma_i}$ represents the identity matrix on $\Gamma_i$ and $\mathcal{I}_{\Omega_j\to\Gamma_i}$ is the interpolation operator between the nodes in $\Omega_j$ and those on $\Gamma_i$ ($i,j=1,2$, $i\not=j$).
System \eqref{eq:systemSchwarz} can be equivalently reformulated only in terms of the interface unknowns $\mathbf{u}_{\Gamma_1}(\bar{\bmu})$ and $\mathbf{u}_{\Gamma_2}(\bar{\bmu})$ as the interface system
\begin{equation}\label{eq:systemSchwarzInterface}
\hspace*{-3mm}
\begin{pmatrix}
{I}_{\Gamma_1} & \mathcal{I}_{\Omega_2\to\Gamma_1} {A}_{\Omega_2}^{-1} {A}_{\Gamma_2} \\[3pt]
\mathcal{I}_{\Omega_1\to\Gamma_2} {A}_{\Omega_1}^{-1} {A}_{\Gamma_1} & {I}_{\Gamma_2}
\end{pmatrix}
\begin{pmatrix}
\mathbf{u}_{\Gamma_1}\!(\bar{\bmu}) \\[3pt]
\mathbf{u}_{\Gamma_2}\!(\bar{\bmu})
\end{pmatrix}
=
\begin{pmatrix}
\mathcal{I}_{\Omega_2\to\Gamma_1} {A}_{\Omega_2}^{-1} \, \mathbf{f}_{\Omega_2}\!(\bar{\bmu}) \\[3pt]
\mathcal{I}_{\Omega_1\to\Gamma_2} {A}_{\Omega_1}^{-1} \, \mathbf{f}_{\Omega_1}\!(\bar{\bmu})
\end{pmatrix} \,  ,
\end{equation}
which can be solved using a suitable matrix-free Krylov method. The expensive part of this iterative algorithm is the solution of the local problems in $\Omega_i$ (corresponding to $A_{\Omega_i}^{-1}$). This can become especially demanding when a new set of parameters $\bmu$ must be considered because, in general, the matrices $A_{\Omega_i}$ and $A_{\Gamma_i}$ can depend on $\bmu$, so that the entire procedure should be re-executed from scratch.

To reduce the computational cost of the coupling procedure \eqref{eq:systemSchwarzInterface}, we exploit the PGD local surrogate models computed in the offline phase. Indeed, we notice that the vector $\mathbf{u}_{\Gamma_i}\!(\bar{\bmu})$ corresponds to the vector of parameters $\bLambda_i$ in \eqref{eq:lambdaRep}, so that the $i$-th equation of the linear system \eqref{eq:systemSchwarzInterface} can be rewritten as
\begin{equation}\label{eq:SchwarzAlgLambda}
{I}_{\Gamma_i}\bLambda_i
=
\mathcal{I}_{\Omega_j\to\Gamma_i}  \, \left( {A}_{\Omega_j}^{-1} \mathbf{f}_{\Omega_j}\!(\bar{\bmu}) 
+  {A}_{\Omega_j}^{-1} (-{A}_{\Gamma_j}\bLambda_j )\, \right)\, .
\end{equation}
Here, ${A}_{\Omega_j}^{-1}\, \mathbf{f}_{\Omega_j}\!(\bar{\bmu})$ corresponds to the vector of the nodal values of the PGD solution of problem~\eqref{eq:sourceProb}, while ${A}_{\Omega_j}^{-1} (-{A}_{\Gamma_j}\bLambda_j)$ corresponds to the linear combination, with coefficient $\bLambda_j$, of the vectors of the nodal values of the PGD solutions of problem~\eqref{eq:boundaryProbs}. All vectors are evaluated for the target parametric values $\bar{\bmu} \in \mathcal{P}$.
Therefore, using the notation introduced in Sect.~\ref{sec:surrogateModels}, we can equivalently write
\begin{equation}\label{eq:SchwarzAlgPGD}
{I}_{\Gamma_i}\bLambda_i
=
\mathcal{I}_{\Omega_j\to\Gamma_i} \, \left( \mathbf{u}_{j,0}^\pgd\,(\bar{\bmu}) 
+ \sum_{q=1}^{N_{\Gamma_j}} \Lambda_j^q \mathbf{u}_{j,q}^\pgd \,(\bar{\bmu}) \, \right)  \, .
\end{equation}
For $j=1,2$, we can define the local PGD operator 
\begin{equation}\label{eq:PGDlambdaOper}
\mathcal{A}_j^\pgd : \bLambda_j \to \sum_{q=1}^{N_{\Gamma_j}} \Lambda_j^q \mathbf{u}_{j,q}^\pgd \,(\bar{\bmu}) \, 
\end{equation}
so that, equation \eqref{eq:SchwarzAlgPGD} and the interface system \eqref{eq:systemSchwarzInterface} can respectively be rewritten in the compact forms
\begin{equation}\label{eq:SchwarzAlgPGD_Op}
{I}_{\Gamma_i}\bLambda_i
=
\mathcal{I}_{\Omega_j\to\Gamma_i}
\left( \mathbf{u}_{j,0}^\pgd \, (\bar{\bmu}) + \mathcal{A}_j^\pgd \, \bLambda_j
\right) \, ,    
\end{equation}
and 
\begin{equation}\label{eq:systemSchwarzInterfacePGD}
\begin{pmatrix}
{I}_{\Gamma_1} & -\mathcal{I}_{\Omega_2\to\Gamma_1} \mathcal{A}_{2}^\pgd \\[3pt]
-\mathcal{I}_{\Omega_1\to\Gamma_2} \mathcal{A}_{1}^\pgd & {I}_{\Gamma_2}
\end{pmatrix}
\begin{pmatrix}
\bLambda_1 \\[3pt]
\bLambda_2
\end{pmatrix}
=
\begin{pmatrix}
\mathcal{I}_{\Omega_2\to\Gamma_1} \mathbf{u}_{2,0}^\pgd \,(\bar{\bmu}) \\[3pt]
\mathcal{I}_{\Omega_1\to\Gamma_2} \mathbf{u}_{1,0}^\pgd\,(\bar{\bmu})
\end{pmatrix} \, .
\end{equation}

Therefore, the online phase of the PGD-enhanced Schwarz method consists in solving \eqref{eq:systemSchwarzInterfacePGD} by an iterative matrix-free method, e.g., GMRES. The computational cost of each GMRES iteration is mainly due to performing the linear combination \eqref{eq:PGDlambdaOper} ($j=1,2$). At convergence, say, at iteration $k=k^*$, the approximation of the solution of the global problem~\eqref{eq:globalProb} for $\bar{\bmu} \in \mathcal{P}$ is given by
\begin{equation}\label{eq:globalSolnA}
	\mathbf{u}^\pgd (\bar{\bmu}) = 
	\begin{cases}
		\mathbf{u}_{1,0}^\pgd (\bar{\bmu}) + \sum_{q=1}^{N_{\Gamma_1}} (\Lambda_1^q)^{(k^*)} \, \mathbf{u}_{1,q}^\pgd (\bar{\bmu}) \quad \text{in } \Omega_1,\\
		\noalign{\vskip5pt}
		\mathbf{u}_{2,0}^\pgd (\bar{\bmu}) + \sum_{q=1}^{N_{\Gamma_2}} (\Lambda_2^q)^{(k^*)} \, \mathbf{u}_{2,q}^\pgd (\bar{\bmu}) \quad \text{in } \Omega_2 \setminus \Omega_{12}.
	\end{cases}
\end{equation}

The approach presented here for the online phase is more efficient than the one based on the clustering of the nodes in the offline phase, see Sect.~\ref{sec:offline} and \cite{deg-cmame-2024}. Indeed, in the latter case, the term $A_{\Omega_j}^{-1}(-A_{\Gamma_j}\bLambda_j)$ corresponds to the vector $\mathbf{u}_{j,\Lambda}^\pgd (\bmu,\bLambda_j)$ of the nodal values of the PGD solution of problem \eqref{eq:boundaryProbsCluster}. Therefore, \eqref{eq:SchwarzAlgPGD} becomes
\begin{equation}\label{eq:SchwarzAlgPGDcluster}
{I}_{\Gamma_i}\bLambda_i
=
\mathcal{I}_{\Omega_j\to\Gamma_i} \, \left( \mathbf{u}_{j,0}^\pgd\,(\bar{\bmu}) 
+ \mathbf{u}_{j,\Lambda}^\pgd (\bmu,\bLambda_j) \, \right)  \, ,
\end{equation}
and the operator $\mathcal{A}_j^\pgd$ must be replaced by $A_{j,\Lambda_j}^\pgd: \, \bLambda_j \to \mathbf{u}_{j,\Lambda}^\pgd (\bar{\bmu},\bLambda_j)$, which returns the nodal values of the PGD model that satisfies problem \eqref{eq:boundaryProbsCluster} for the boundary parameters $\bLambda_j$ and for $\bar{\bmu}$.
Since the values $\bLambda_j^{(k)}$ computed at the $k$-th GMRES iteration do not generally coincide with those obtained when discretizing the parametric space $\mathcal{Q}_j$, applying the local operator $A_{j,\Lambda_j}^\pgd$ requires linearly interpolating the parametric modes that depend on $\bLambda_j$ and that are associated with the available values closest to $\bLambda_j^{(k)}$. Although this operation is computationally inexpensive, its repeated execution involves a higher computational cost than the linear combination \eqref{eq:PGDlambdaOper}.

Moreover, for the approach based on clustering the interface nodes, the parametric domain $\mathcal{Q}_j$ is fixed during the offline phase and, in principle, it could happen that some of the values $\bLambda_j^{(k)}$ computed by GMRES fall outside of the parametric space $\mathcal{Q}_j$, thus jeopardizing the convergence of the iterative scheme. This cannot occur using the approach described in Sect.~\ref{sec:offline}. Indeed, it is always possible to compute the linear combination \eqref{eq:PGDlambdaOper} for any given coefficients $\bLambda_j^{(k)}$ once the \textit{PGD basis functions} $\mathbf{u}_{j,q}^\pgd(\bar{\bmu})$ are available. The resulting online coupling algorithm is thus computationally less expensive and more robust.

\section{Numerical results}\label{sec:numericalResults}

In this section, we reproduce the tests studied in Sects. 5.1 and 5.3 of \cite{deg-cmame-2024} using the approach explained in this paper to assess its performance versus the technique based on clustering the boundary nodes. Due to the restrictions on the number of pages, we outline the definition of the problems and refer to \cite{deg-cmame-2024} for a detailed description.
For all tests, in the offline phase, the local parametric problems \eqref{eq:sourceProb} and \eqref{eq:boundaryProbs} are solved using the encapsulated PGD toolbox \cite{epgd-2020} with tolerance $10^{-4}$ to stop the PGD enrichment process, and a compression algorithm \cite{mzh-cmame-2015} with tolerance $10^{-3}$ is applied to eliminate redundant modes. In the online phase, the interface system \eqref{eq:systemSchwarzInterfacePGD} is solved by GMRES with stopping tolerance $10^{-6}$ on the relative residual.

\subsection{Bidomain diffusion problem with analytical solution}

Consider problem \eqref{eq:globalProb} with $\nu(\mu) = 1+\mu x$, $\bx=(x,y)\in\Omega$ and $\mu \in \mathcal{P}:=[1,50]$ being a scalar parameter. The domain is $\Omega = (0,2)\times (0,1)$, split into the two overlapping subdomains $\Omega_1 = (0, 1.05) \times (0, 1)$ and $\Omega_2 = (0.95, 2) \times (0, 1)$, and the source term $s(\mu)$ is chosen such that the exact solution is $u(\mu) = \sin(2\pi x) \, \sin(2\pi y) + \frac{\mu}{2} \, xy(y-1)(x-2)$. The parametric interval $\mathcal{P}$ is discretized using step $h_\mu = 10^{-3}$, while Lagrangian $\mathbb{Q}_1$ elements are used for the space discretization with mesh size $h=0.05$. The local meshes coincide in the overlap that has width $2h$. Considering this space discretization, there are 19 nodes on each interface so that, following the procedure in Sect.~\ref{sec:offline}, one must solve $19$ problems like \eqref{eq:boundaryProbs} and one problem \eqref{eq:sourceProb} in each subdomain. The total number of PGD modes computed in $\Omega_1$ is 68 (106 before compression) and 56 (62 before compression) in $\Omega_2$, with total CPU time of the offline phase of 6.56~s (averaged over 10 runs).
In the online phase for the case $\mu=3$, GMRES converges in 9 iterations with computational time of about $2 \times 10^{-2}$~s. To assess the accuracy of the DD-PGD approach, we compute the relative error in $L^2$ norm between the exact solution and the DD-PGD solution for selected values of the parameter $\mu$ and we compare them to the errors for the full-order finite element solution $u^h_\Omega$. The results reported in Table \ref{tab:diffusion} show that the accuracy provided by the two methods is almost identical.

\begin{table}
\begin{center}
\begin{tabular}{ccc}
$\mu$ & $\|u^\pgd-u\|_{L^2(\Omega)}/\|u\|_{L^2(\Omega)}$ & $\|u_\Omega^h-u\|_{L^2(\Omega)}/\|u\|_{L^2(\Omega)}$ \\
\hline
 3 & $9.08 \times 10^{-3}$ & $9.07 \times 10^{-3}$ \\
30 & $3.27 \times 10^{-3}$ & $3.27 \times 10^{-3}$ \\
\hline
\end{tabular}
\end{center}
\caption{Relative error in $L^2$ norm between the exact solution $u$, the DD-PGD solution $u^\pgd$ and the full-order finite element solution $u_\Omega^h$.}
\label{tab:diffusion}
\end{table}

Comparing these results to those in Sect. 5.1 of \cite{deg-cmame-2024}, we remark that the offline computational cost is now reduced by approximately 17 times without loss of accuracy in the online phase. This is due to the fact that many unneeded parametric combinations occurring in \eqref{eq:boundaryProbsCluster} are avoided with the present approach. The number of modes is also lower, thus reducing memory storage. Meanwhile, the online phase is twice as fast and requires fewer iterations with both advantages being attributable to the absence of an interpolation procedure for the boundary parameters.

\subsection{Multi-domain diffusion problem}

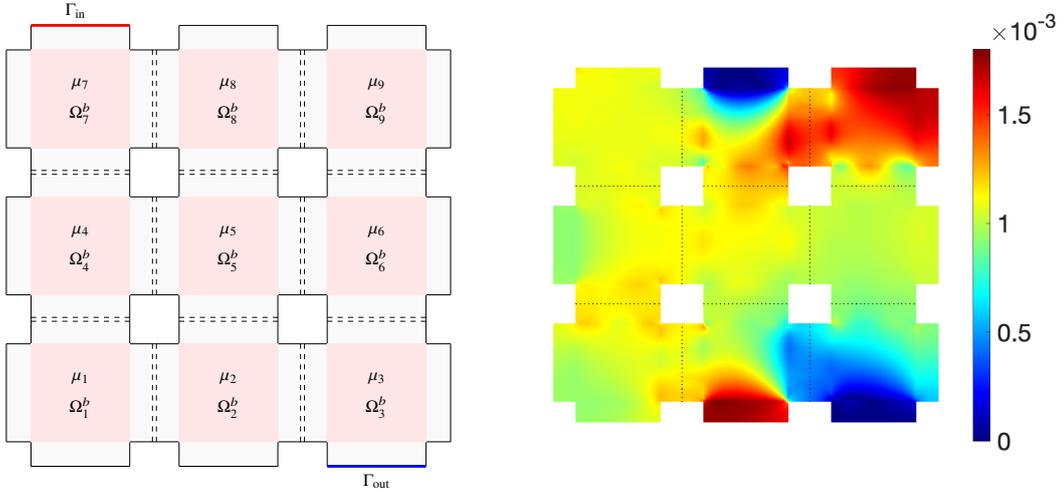
\begin{figure}[thb]
\begin{center}
\resizebox{0.4\textwidth}{!}{
\begin{tikzpicture}
        \def\s{2}
        \begin{scope}
            \pateraWing{1}{0}{1}{0}
            \pateraBlock{1}{1}{1}{2}{2}{1}{1}
        \end{scope}

        \begin{scope}[shift={(1.5*\s, 0)}]
            \pateraWing{1}{0}{0}{0}
            \pateraBlock{2}{1}{2}{2}{2}{1}{1}
        \end{scope}

        \begin{scope}[shift={(3*\s, 0)}]
            \pateraWing{1}{0}{0}{1}
            \pateraBlock{3}{1}{2}{1}{2}{1}{1}
            \draw[ultra thick, blue] (0, -0.25*\s) -- (1*\s, -0.25*\s);
        \end{scope}

        \begin{scope}[shift={(0, 1.5*\s)}]
            \pateraWing{0}{0}{1}{0}
            \pateraBlock{4}{2}{1}{2}{2}{1}{1}
        \end{scope}

        \begin{scope}[shift={(1.5*\s, 1.5*\s)}]
            \pateraWing{0}{0}{0}{0}
            \pateraBlock{5}{2}{2}{2}{2}{1}{1}
        \end{scope}

        \begin{scope}[shift={(3*\s, 1.5*\s)}]
            \pateraWing{0}{0}{0}{1}
            \pateraBlock{6}{2}{2}{1}{2}{1}{1}
        \end{scope}

        \begin{scope}[shift={(0, 3*\s)}]
            \pateraWing{0}{1}{1}{0}
            \pateraBlock{7}{2}{1}{2}{1}{1}{1}
            \draw[ultra thick, red] (0, 1.25*\s) -- (1*\s, 1.25*\s);
        \end{scope}

        \begin{scope}[shift={(1.5*\s, 3*\s)}]
            \pateraWing{0}{1}{0}{0}
            \pateraBlock{8}{2}{2}{2}{1}{1}{1}
        \end{scope}

        \begin{scope}[shift={(3*\s, 3*\s)}]
            \pateraWing{0}{1}{0}{1}
            \pateraBlock{9}{2}{2}{1}{1}{1}{1}
        \end{scope}

        \node at (0.45*\s, 4.4*\s) {$\Gamma_{\text{in}}$};

        \node at (3.5*\s, -0.4*\s) {$\Gamma_{\text{out}}$};

    \end{tikzpicture}
    }
\includegraphics[width=0.53\textwidth]{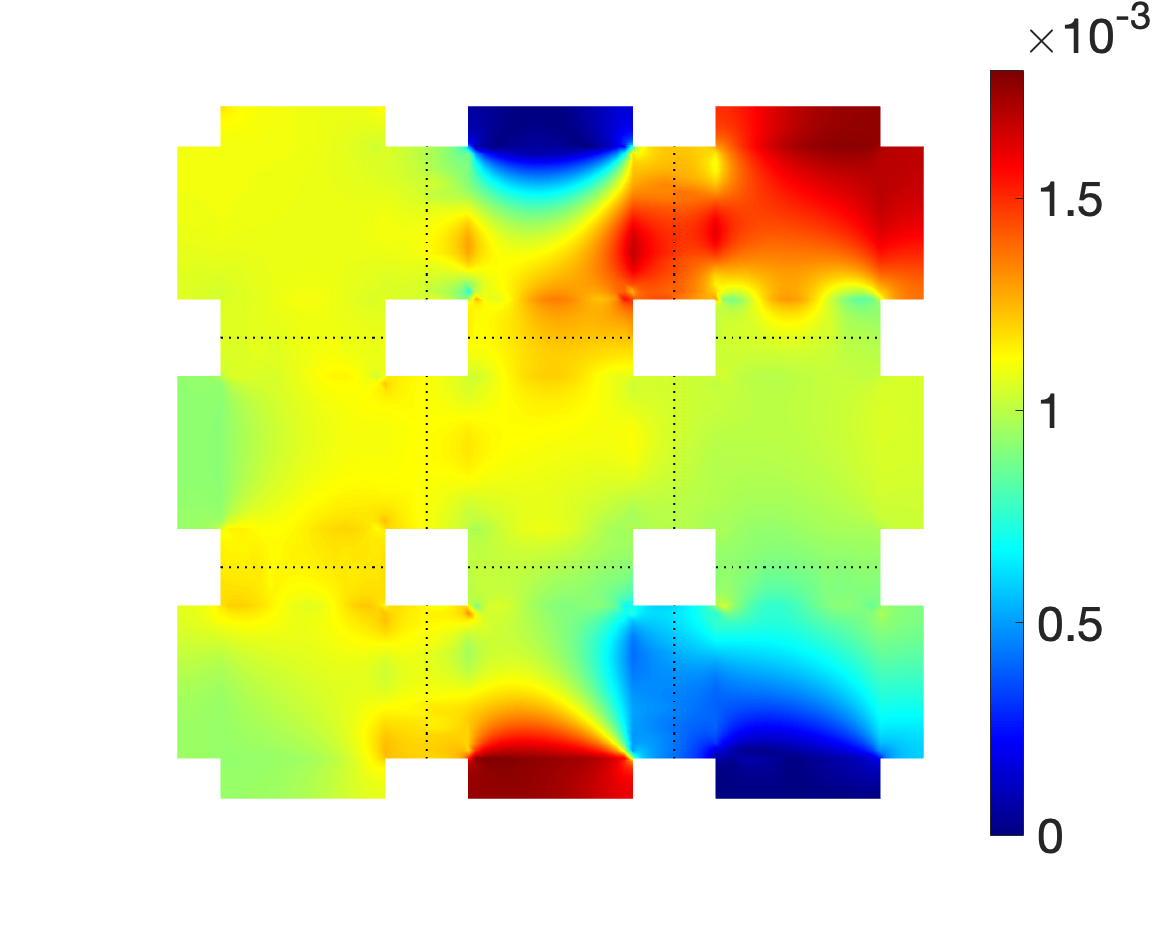}
\end{center}
\caption{Computational domain (left) and scaled error $|u^\pgd(\bmu)-u^h(\mu)| / \max_\Omega | u^h(\bmu)|$ for the first set of parameters reported in Table~\ref{tab:testPatera}.}
\label{fig:testPatera}
\end{figure}

We consider the diffusion equation \eqref{eq:globalProb}$_1$ on the domain $\Omega$ shown in Fig.~\ref{fig:testPatera} (left) with boundary conditions $u(\bmu)=0$ on $\Gamma_{\text{out}}$, $\nu(\bmu)\nabla u(\bmu) \cdot\bn = 1$ on $\Gamma_{\text{in}}$ and $\nabla u(\bmu) \cdot \bn = 0$ on $\partial\Omega \setminus (\Gamma_{\text{in}} \cup \Gamma_{\text{out}})$. In this test case, proposed in \cite{ep-ijnme-2013}, each subdomain $\Omega_i^b$ is characterized by a parametric diffusion coefficient $\nu(\bmu) = \mu_i$ ($i=1,\ldots,9$), while $\nu(\bmu)=1$ in the rest of the domain, so that there are 9 scalar parameters $\mu_i \in \mathcal{P} := [5\times 10^{-2},10]$. As we remarked in \cite{deg-cmame-2024}, the computational domain can be reconstructed by composing four reference subdomains $\hat{\Omega}_i$ ($i=1,\ldots,4$) upon suitable rotation and/or translation. Therefore, in the offline phase, the local surrogate models are computed only in those reference domains considering discretization step $h_\mu = 10^{-3}$ for $\mathcal{P}$, and $\mathbb{Q}_1$ Lagrangian elements for the finite element discretization with mesh size $h=0.0125$ along the external boundary of each subdomain and $h = 0.05$ within the bulk regions $\Omega_i^b$. This space discretization leads to 42 boundary problems \eqref{eq:boundaryProbs} in $\hat{\Omega}_1$ and in $\hat{\Omega}_4$, 63 problems in $\hat{\Omega}_2$, and 84 in $\hat{\Omega}_3$. In $\hat{\Omega}_4$, an additional problem \eqref{eq:sourceProb} is considered to account for the Neumann boundary condition on $\Gamma_{\text{in}}$. The number of PGD modes computed in the offline phase in $\hat{\Omega}_i$ ($i=1,\ldots,4$) is, respectively, 190 (312 before compression), 272 (469), 336 (604) and 195 (321), with a total computational cost of approximately 35~s.
These are significantly fewer modes than those obtained by the clustering procedure that resulted, respectively, in 236 (370 before compression), 386 (672), 509 (930) and 244 (419) modes with computational cost approximately 25 times higher.
The number of GMRES iterations and the CPU time obtained for two sets of parameters $\bmu$ in the online phase are reported in Table~\ref{tab:testPatera}, together with the relative errors in $l^\infty$ norm versus a full-order finite element solution $u^h_\Omega(\bmu)$. The error distribution for the first set of parameters is also shown in Fig.~\ref{fig:testPatera} (right).

\begin{table}[hbt]
\begin{center}
\begin{tabular}{ccccccccccccc}
 & $\mu_1$ & $\mu_2$ & $\mu_3$ & $\mu_4$ & $\mu_5$ & $\mu_6$ & $\mu_7$ & $\mu_8$ & $\mu_9$ & \# iter. & CPU time & $\frac{\|u^\pgd(\boldsymbol{\mu}) - u^h(\boldsymbol{\mu})\|_{l^\infty(\Omega)}}{\|u^h(\boldsymbol{\mu})\|_{l^\infty(\Omega)}}$\\
\cline{2-13}
Test 1: \quad & 0.1 & 0.2 & 0.4 & 0.8 & 1.6 & 3.2 & 6.4 & 0.1 & 0.2 & 93 &  1.26~s & $1.8 \times 10^{-3}$ \\
Test 2: \quad & 4.9 & 4.7 & 4.8 & 5.2 & 5.0 & 4.9 & 5.5 & 5.3 & 5.1 & 57 & 0.79~s & $3.4 \times 10^{-4}$ \\
\cline{2-13}
\end{tabular}
\end{center}
\caption{Number of GMRES iterations, CPU time and relative errors obtained in the online phase for two sets of parametric values.}
\label{tab:testPatera}
\end{table}

We remark that the number of GMRES iterations is significantly lower than in the case of clustering of the nodes where 302 and 126 iterations were needed to converge in Test 1 and 2, respectively. The CPU time is also reduced by approximately 15 and 10 times.
It is also worth noticing that the number of performed GMRES iterations (93 and 57, respectively) is now comparable to the 95 and 56 iterations required by the overlapping Schwarz algorithm coupled with the high-fidelity finite element solver, while the proposed strategy offers a wall-clock speed-up of approximately 130 times.
The fewer iterations are likely due to the avoidance of the interpolation procedure that was needed by the previous clustering approach.

\section{Conclusions}

In this work, we presented an improved algorithm that combines overlapping DD and physics-based PGD MOR to solve parametric linear elliptic problems. The proposed approach exploits the basis functions that span the space of the traces of the finite element solution on the domain interfaces to construct a \textit{PGD basis} for the local surrogate models, with a significant reduction of the dimensionality of the parametric subproblems to be solved. These are then efficiently added through parametric linear combinations to obtain the local surrogate models needed at each iteration of the Schwarz algorithm in the online phase. The resulting DD-PGD algorithm outperforms the approach in \cite{deg-cmame-2024} in terms of simplicity of implementation and of computational cost both in the online and in the offline phase.

\bigskip

\textbf{Acknowledgements.}
MD acknowledges EPSRC grant EP/V027603/1, BJE the EPSRC Doctoral Training Partnership grant EP/W523987/1, and MG the Spanish Ministry of Science, Innovation and Universities and State Research Agency MICIU/AEI/10.13039/501100011033 (Grants PID2020-113463RB-C33, CEX2018-000797-S) and Generalitat de Catalunya (Grant 2021-SGR-01049). MG is Fellow of the Serra H\'unter Programme of the Generalitat de Catalunya.

\smallskip

\textbf{Competing Interests} {The authors have no conflicts of interest to declare that are relevant to the content of this paper.}


\begin{thebibliography}{99.}
%
\bibitem{buhr-2021}
Buhr, A., Iapichino, L., Ohlberger, M., Rave, S., Schindler, F., Smetana, K.: Localized model reduction for parameterized problems. In: Benner, P., et al. (eds.) Snapshot-Based Methods and Algorithms (Volume 2), pp. 245--306, De Gruyter (2021).
%
\bibitem{ckl-pgdBook}
Chinesta, F., Keunings, R., Leygue, A.: The Proper Generalized Decomposition for Advanced Numerical Simulations. Springer (2014).
%
\bibitem{epgd-2020}
D{\'\i}ez, P., Zlotnik, S., Garc{\`\i}a-Gonz\'alez, A., Huerta, A.: Encapsulated {PGD} algebraic toolbox operating with high-dimensional data. Arch. Comput. Methods Eng. \textbf{27}, 1321--1336 (2020).
%
\bibitem{deg-cmame-2024}
Discacciati, M., Evans, B.J., Giacomini, M.: An overlapping domain decomposition method for the solution of parametric elliptic problems via proper generalized decomposition. Comput. Methods Appl. Mech. Engrg. \textbf{418}, 116484 (2024).
%
\bibitem{deg-stokesDarcy-2024}
Discacciati, M., Evans, B.J., Giacomini, M.: An overlapping domain decomposition method for the Stokes and Stokes-Darcy problems via PGD. In preparation (2024).
%
\bibitem{ep-ijnme-2013}
Erftang, J.L., Patera, A.T.: Port reduction in parametrized component static condensation and approximation and a posteriori error estimation. Internat. J. Numer. Methods Engrg. \textbf{96}(5), 269--302 (2013).
%
\bibitem{mzh-cmame-2015}
Modesto, D., Zlotnik, S., Huerta, A.: Proper generalized decomposition for parametrized {H}elmholtz problems in heterogeneous and unbounded domains: Application to harbor agitation. Comput. Methods Appl. Mech. Engrg. \textbf{295}, 127--149 (2015).
%
\bibitem{smith-ddBook}
Smith, B., Bj{\o}rstad, P., Gropp, W.: Domain Decomposition: Parallel Multilevel Methods for Elliptic Partial Differential Equations. Cambridge University Press (1996).
%
\end{thebibliography}
\end{document}